\documentclass[12pt]{article}

\usepackage{amsmath,amssymb,amsthm,amscd,a4wide,upref}

\usepackage[enableskew]{youngtab}

\usepackage{hyperref}
\hypersetup{colorlinks,citecolor=blue,filecolor=black,linkcolor=blue,urlcolor=blue}
\usepackage{enumerate}
\usepackage{mathtools}
\usepackage{enumitem}

\usepackage{lmodern}     
\usepackage[T1]{fontenc}

\begin{document}


\newcommand{\ad}{{\rm ad}}
\newcommand{\cri}{{\rm cri}}
\newcommand{\row}{{\rm row}}
\newcommand{\col}{{\rm col}}
\newcommand{\Ann}{{\rm{Ann}\ts}}
\newcommand{\End}{{\rm{End}\ts}}
\newcommand{\Rep}{{\rm{Rep}\ts}}
\newcommand{\Hom}{{\rm{Hom}}}
\newcommand{\Mat}{{\rm{Mat}}}
\newcommand{\ch}{{\rm{ch}\ts}}
\newcommand{\chara}{{\rm{char}\ts}}
\newcommand{\diag}{{\rm diag}}
\newcommand{\st}{{\rm st}}
\newcommand{\non}{\nonumber}
\newcommand{\wt}{\widetilde}
\newcommand{\wh}{\widehat}
\newcommand{\ol}{\overline}
\newcommand{\ot}{\otimes}
\newcommand{\la}{\lambda}
\newcommand{\La}{\Lambda}
\newcommand{\De}{\Delta}
\newcommand{\al}{\alpha}
\newcommand{\be}{\beta}
\newcommand{\ga}{\gamma}
\newcommand{\Ga}{\Gamma}
\newcommand{\ep}{\epsilon}
\newcommand{\ka}{\kappa}
\newcommand{\vk}{\varkappa}
\newcommand{\si}{\sigma}
\newcommand{\vs}{\varsigma}
\newcommand{\vp}{\varphi}
\newcommand{\ta}{\theta}
\newcommand{\de}{\delta}
\newcommand{\ze}{\zeta}
\newcommand{\om}{\omega}
\newcommand{\Om}{\Omega}
\newcommand{\ee}{\epsilon^{}}
\newcommand{\su}{s^{}}
\newcommand{\hra}{\hookrightarrow}
\newcommand{\rar}{\rightarrow}
\newcommand{\lar}{\leftarrow}
\newcommand{\ve}{\varepsilon}
\newcommand{\pr}{^{\tss\prime}}
\newcommand{\ts}{\,}
\newcommand{\vac}{\mathbf{1}}
\newcommand{\vacu}{|0\rangle}
\newcommand{\di}{\partial}
\newcommand{\qin}{q^{-1}}
\newcommand{\tss}{\hspace{1pt}}
\newcommand{\Sr}{ {\rm S}}
\newcommand{\U}{ {\rm U}}
\newcommand{\BL}{ {\overline L}}
\newcommand{\BE}{ {\overline E}}
\newcommand{\BP}{ {\overline P}}
\newcommand{\AAb}{\mathbb{A}\tss}
\newcommand{\CC}{\mathbb{C}\tss}
\newcommand{\KK}{\mathbb{K}\tss}
\newcommand{\QQ}{\mathbb{Q}\tss}
\newcommand{\SSb}{\mathbb{S}\tss}
\newcommand{\TT}{\mathbb{T}\tss}
\newcommand{\ZZ}{\mathbb{Z}\tss}
\newcommand{\DY}{ {\rm DY}}
\newcommand{\X}{ {\rm X}}
\newcommand{\Y}{ {\rm Y}}
\newcommand{\Z}{{\rm Z}}
\newcommand{\ZX}{{\rm ZX}}
\newcommand{\Ac}{\mathcal{A}}
\newcommand{\Lc}{\mathcal{L}}
\newcommand{\Mc}{\mathcal{M}}
\newcommand{\Pc}{\mathcal{P}}
\newcommand{\Qc}{\mathcal{Q}}
\newcommand{\Rc}{\mathcal{R}}
\newcommand{\Sc}{\mathcal{S}}
\newcommand{\Tc}{\mathcal{T}}
\newcommand{\Bc}{\mathcal{B}}
\newcommand{\Ec}{\mathcal{E}}
\newcommand{\Fc}{\mathcal{F}}
\newcommand{\Gc}{\mathcal{G}}
\newcommand{\Hc}{\mathcal{H}}
\newcommand{\Uc}{\mathcal{U}}
\newcommand{\Vc}{\mathcal{V}}
\newcommand{\Wc}{\mathcal{W}}
\newcommand{\Yc}{\mathcal{Y}}
\newcommand{\Cl}{\mathcal{C}l}
\newcommand{\Ar}{{\rm A}}
\newcommand{\Br}{{\rm B}}
\newcommand{\Ir}{{\rm I}}
\newcommand{\Fr}{{\rm F}}
\newcommand{\Jr}{{\rm J}}
\newcommand{\Or}{{\rm O}}
\newcommand{\GL}{{\rm GL}}
\newcommand{\Spr}{{\rm Sp}}
\newcommand{\Rr}{{\rm R}}
\newcommand{\Zr}{{\rm Z}}
\newcommand{\gl}{\mathfrak{gl}}
\newcommand{\middd}{{\rm mid}}
\newcommand{\ev}{{\rm ev}}
\newcommand{\Pf}{{\rm Pf}}
\newcommand{\Norm}{{\rm Norm\tss}}
\newcommand{\oa}{\mathfrak{o}}
\newcommand{\spa}{\mathfrak{sp}}
\newcommand{\osp}{\mathfrak{osp}}
\newcommand{\f}{\mathfrak{f}}
\newcommand{\g}{\mathfrak{g}}
\newcommand{\h}{\mathfrak h}
\newcommand{\n}{\mathfrak n}
\newcommand{\m}{\mathfrak m}
\newcommand{\se}{\mathfrak{s}}
\newcommand{\z}{\mathfrak{z}}
\newcommand{\Zgot}{\mathfrak{Z}}
\newcommand{\p}{\mathfrak{p}}
\newcommand{\sll}{\mathfrak{sl}}
\newcommand{\agot}{\mathfrak{a}}
\newcommand{\bgot}{\mathfrak{b}}
\newcommand{\qdet}{ {\rm qdet}\ts}
\newcommand{\Ber}{ {\rm Ber}\ts}
\newcommand{\HC}{ {\mathcal HC}}
\newcommand{\cdet}{{\rm cdet}}
\newcommand{\rdet}{{\rm rdet}}
\newcommand{\tr}{ {\rm tr}}
\newcommand{\gr}{ {\rm gr}\ts}
\newcommand{\str}{ {\rm str}}
\newcommand{\loc}{{\rm loc}}
\newcommand{\Gr}{{\rm G}}
\newcommand{\sgn}{ {\rm sgn}\ts}
\newcommand{\sign}{{\rm sgn}}
\newcommand{\ba}{\bar{a}}
\newcommand{\bb}{\bar{b}}
\newcommand{\bi}{\bar{\imath}}
\newcommand{\bj}{\bar{\jmath}}
\newcommand{\bk}{\bar{k}}
\newcommand{\bl}{\bar{l}}
\newcommand{\bp}{\bar{p}}
\newcommand{\hb}{\mathbf{h}}
\newcommand{\Sym}{\mathfrak S}
\newcommand{\fand}{\quad\text{and}\quad}
\newcommand{\Fand}{\qquad\text{and}\qquad}
\newcommand{\For}{\qquad\text{or}\qquad}
\newcommand{\for}{\quad\text{or}\quad}
\newcommand{\grpr}{{\rm gr}^{\tss\prime}\ts}
\newcommand{\degpr}{{\rm deg}^{\tss\prime}\tss}
\newcommand{\bideg}{{\rm bideg}\ts}

\renewcommand{\theequation}{\arabic{section}.\arabic{equation}}

\numberwithin{equation}{section}

\newtheorem{thm}{Theorem}[section]
\newtheorem{lem}[thm]{Lemma}
\newtheorem{prop}[thm]{Proposition}
\newtheorem{cor}[thm]{Corollary}
\newtheorem{conj}[thm]{Conjecture}
\newtheorem*{mthm}{Main Theorem}
\newtheorem*{mthma}{Theorem A}
\newtheorem*{mthmb}{Theorem B}
\newtheorem*{mthmc}{Theorem C}
\newtheorem*{mthmd}{Theorem D}

\theoremstyle{definition}
\newtheorem{defin}[thm]{Definition}

\theoremstyle{remark}
\newtheorem{remark}[thm]{Remark}
\newtheorem{example}[thm]{Example}
\newtheorem{examples}[thm]{Examples}

\newcommand{\bth}{\begin{thm}}
\renewcommand{\eth}{\end{thm}}
\newcommand{\bpr}{\begin{prop}}
\newcommand{\epr}{\end{prop}}
\newcommand{\ble}{\begin{lem}}
\newcommand{\ele}{\end{lem}}
\newcommand{\bco}{\begin{cor}}
\newcommand{\eco}{\end{cor}}
\newcommand{\bde}{\begin{defin}}
\newcommand{\ede}{\end{defin}}
\newcommand{\bex}{\begin{example}}
\newcommand{\eex}{\end{example}}
\newcommand{\bes}{\begin{examples}}
\newcommand{\ees}{\end{examples}}
\newcommand{\bre}{\begin{remark}}
\newcommand{\ere}{\end{remark}}
\newcommand{\bcj}{\begin{conj}}
\newcommand{\ecj}{\end{conj}}

\newcommand{\bal}{\begin{aligned}}
\newcommand{\eal}{\end{aligned}}
\newcommand{\beq}{\begin{equation}}
\newcommand{\eeq}{\end{equation}}
\newcommand{\ben}{\begin{equation*}}
\newcommand{\een}{\end{equation*}}

\newcommand{\bpf}{\begin{proof}}
\newcommand{\epf}{\end{proof}}

\def\beql#1{\begin{equation}\label{#1}}

\newcommand{\Res}{\mathop{\mathrm{Res}}}

\title{\Large\bf Representations of the super-Yangian of type $B(n,m)$}

\author{Alexander Molev\ \   and\ \   Eric Ragoucy}

\date{} 
\maketitle


\begin{abstract}
We are concerned with finite-dimensional irreducible representations of
the Yangians associated with the orthosymplectic Lie superalgebras $\osp_{2n+1|2m}$.
Every such representation is highest weight and we use embedding theorems and
odd reflections of Yangian type to derive necessary
conditions for an irreducible highest weight representation to be finite-dimensional.
We conjecture that these conditions are also sufficient. We prove the conjecture
in the case $n=1$ and arbitrary $m\geqslant 1$.
\end{abstract}



%

\section{Introduction}\label{sec:int}
\setcounter{equation}{0}

We continue the investigation of finite-dimensional irreducible representations of
the orthosymplectic Yangians initiated in \cite{m:ry}
and \cite{m:rs}, where classification theorems were proved in the particular cases
of $\osp_{1|2m}$ and $\osp_{2|2m}$.
We consider the Yangians associated with the Lie superalgebras $\osp_{2n+1|2m}$
which form the series $B(n,m)$ in the Kac classification of simple Lie superalgebras.

We will work with the extended Yangian $\X(\osp_{2n+1|2m})$ and use its $R$-matrix
definition as given by Arnaudon {\it et al.\/}~\cite{aacfr:rp}.
By a standard argument, every finite-dimensional irreducible representation of
$\X(\osp_{2n+1|2m})$ is a highest weight representation. It is isomorphic
to the irreducible quotient $L(\la(u))$ of the Verma module $M(\la(u))$
associated with a tuple
\beql{lau}
\la(u)=\big(\la_1(u),\dots,\la_{m}(u),\la_{m+1}(u),\dots,\la_{m+n+1}(u)\big)
\eeq
of formal series
$\la_i(u)\in 1+u^{-1}\CC[[u^{-1}]]$. The tuple is called the {\em highest weight}
of the representation.
Both the Yangian $\Y(\gl_{n|m})$
associated with the general linear Lie superalgebra
$\gl_{n|m}$ and the extended Yangian $\X(\oa_{2n+1})$ associated
with the orthogonal Lie algebra $\oa_{2n+1}$ can be regarded as subalgebras
of the orthosymplectic Yangian via embeddings
\beql{embedya}
\Y(\gl_{n|m})\hra \X(\osp_{2n+1|2m})\Fand \X(\oa_{2n+1})\hra\X(\osp_{2n+1|2m}).
\eeq
Hence, by considering the cyclic spans of the highest vector of $L(\la(u))$
with respect to these subalgebras, and using the finite-dimensionality conditions
for highest weight representations of $\Y(\gl_{n|m})$ from \cite{zh:sy}
and of $\X(\oa_{2n+1})$ from \cite{amr:rp}, we get certain necessary conditions
for the representation $L(\la(u))$ to be finite-dimensional.

Moreover, an additional necessary condition is obtained by using
a sequence of {\em odd reflections} as introduced in \cite{m:or}; see also \cite{l:no}.
These are transformations which apply to
pairs $(\al(u),\be(u))$ of formal series of the form
\ben
\bal
\al(u)&=(1+\al_1u^{-1})\dots (1+\al_pu^{-1})\ts\ga(u),\\
\be(u)&=(1+\be_1u^{-1})\dots (1+\be_p\tss u^{-1})\ts\ga(u),
\eal
\een
where $\al_i\ne\be_j$ for all $i,j$, and $\ga(u)\in 1+u^{-1}\CC[[u^{-1}]]$.
The odd reflection is the transformation
\beql{oddrefl}
\big(\al(u),\be(u)\big)\mapsto \big(\be^{[1]}(u),\al^{[1]}(u)\big),
\eeq
where
\ben
\bal
\al^{[1]}(u)&=\big(1+(\al_1+1)\tss u^{-1}\big)\dots \big(1+(\al_p+1)\tss u^{-1}\big)\ts\ga(u),\\
\be^{[1]}(u)&=\big(1+(\be_1+1)\tss u^{-1}\big)\dots \big(1+(\be_p+1)\tss u^{-1}\big)\ts\ga(u).
\eal
\een

Now the additional series $\la^{[n]}_{m}(u)$ for the tuple \eqref{lau}
is obtained by applying the sequence of odd reflections \eqref{oddrefl}
beginning with
\beql{lammpo}
\big(\la_m(u),\la_{m+1}(u)\big)\mapsto \big(\la^{[1]}_{m+1}(u),\la^{[1]}_{m}(u)\big),
\eeq
then using the values of $i=1,\dots,n-1$ and applying \eqref{oddrefl}
consecutively
to determine new series $\la^{[2]}_{m}(u),\dots,\la^{[n]}_{m}(u)$ by
\beql{oddlami}
\big(\la^{[i]}_{m}(u),\la_{m+i+1}(u)\big)\mapsto \big(\la^{[1]}_{m+i+1}(u),\la^{[i+1]}_{m}(u)\big).
\eeq
Note that the required form of the series involved in the transformations
is guaranteed by the necessary conditions pointed out above.
We may now use the embedding
\beql{embospyang}
\X(\osp_{1|2})\hra \X(\osp_{2n+1|2m})
\eeq
and the results of \cite{m:ry} to derive the following.

\bth\label{thm:necess}
If the representation $L(\la(u))$ of the Yangian $\X(\osp_{2n+1|2m})$
is finite-dimensional then
\beql{condgl}
\frac{\la_{i+1}(u)}{\la_{i}(u)}=\frac{P_i(u+1)}{P_i(u)},\qquad\text{for}\quad i=1,\dots,m-1,\qquad
\eeq
\beql{condo}
\frac{\la_{j}(u)}{\la_{j+1}(u)}=\frac{P_j(u+\si_j)}{P_j(u)},\qquad\text{for}\quad j=m+1,\dots,m+n,
\eeq
with $\si_j=1$ for $j<m+n$ and $\si_{m+n}=1/2$, and
\beql{condosp}
\frac{\la_{m+n+1}(u)}{\la^{[n]}_{m}(u)}=\frac{P_m(u+1)}{P_m(u)},
\eeq
for some monic polynomials $P_1(u),\dots,P_{m+n}(u)$ in $u$.
\eth

We will give a detailed proof of Theorem~\ref{thm:necess} in Sec.~\ref{subsec:thmnec}.
Then we bring some evidence provided by particular cases to support the following
conjecture.

\bcj\label{conj:fdim}
The conditions on the highest weight $\la(u)$ stated in Theorem~\ref{thm:necess}
are sufficient for the representation $L(\la(u))$ to be finite-dimensional.
\ecj

Conjecture~\ref{conj:fdim} holds for $n=0$ due to the results of \cite{m:ry},
where conditions \eqref{condo} and \eqref{condosp} do not occur.
We will prove the following theorem in Sec.~\ref{sec:ct}.

\bth\label{thm:none}
Conjecture~\ref{conj:fdim} holds for $n=1$ and arbitrary $m\geqslant 1$.
\eth

Before specializing to the case $n=1$, we give necessary and
sufficient conditions for the representations $L(\la(u))$ of $\X(\osp_{2n+1|2m})$
with {\em linear highest weights}
\beql{linhw}
\la_i(u)=1+\la_i\ts u^{-1},\qquad i=1,\dots,m+n+1,\quad\la_i\in\CC,
\eeq
to be finite-dimensional (Theorem~\ref{thm:linear}). They will show that Conjecture~\ref{conj:fdim} holds
in this case.
We then use the linear highest weight modules to prove Conjecture~\ref{conj:fdim}
in the case of {\em generic} highest weights (Corollary~\ref{cor:genhw}).

A key role in the proof of Theorem~\ref{thm:linear} will be played by
the {\em polynomial representations} of $\X(\osp_{2n+1|2m})$ which are
produced in a way similar to the Yangian $\Y(\gl_{n|m})$, as reviewed e.g. in
\cite[Appendix]{m:rs}.
Moreover, the proof will also use a special
fundamental representation of the Lie superalgebra $\osp_{2n+1|2m}$
which admits an extension to the Yangian, as shown in Lemma~\ref{lem:fundosc}.

\medskip

This work was completed during the first author's visit to
the {\it Laboratoire d'Annecy-le-Vieux de Physique Th\'eorique\/}.
He thanks the lab for the support and hospitality.

\section{Definitions and preliminaries}
\label{sec:ns}

For given positive integers $m$ and $n$ consider
the {\em parity sequences} $\se=\se_1\dots\se_{m+n}$ of length $m+n$,
where each term $\se_i$ is $0$ or $1$, and the total number of zeros is $n$.
For notational convenience we also set $\se_{m+n+1}=0$ and regard this as
an additional entry of $\se$ which will not vary.
The {\em standard sequence} $\se^{\st}=1\dots 1\tss 0\dots 0$
is defined by $\se_i=1$ for $i=1,\dots,m$
and $\se_i=0$ for $i=m+1,\dots,m+n$.

Suppose a parity sequence $\se$ is fixed. We will simply write $\bi$ to
denote its $i$-th term $\se_i$. Introduce the
involution $i\mapsto i\pr=2n+2m-i+2$ on
the set $\{1,2,\dots,2n+2m+1\}$ and set $\overline{i'}=\bi$ for $i=1,\dots,m+n+1$.
Consider the $\ZZ_2$-graded vector space $\CC^{2n+1|2m}$ over $\CC$ with the basis
$e_1,e_2,\dots,e_{1'}$, where the parity of the basis vector
$e_i$ is defined to be $\bi\mod 2$.
Accordingly, equip
the endomorphism algebra $\End\CC^{2n+1|2m}$ with the $\ZZ_2$-gradation, where
the parity of the matrix unit $e_{ij}$ is found by
$\bi+\bj\mod 2$.

We will consider even square matrices with entries in $\ZZ_2$-graded algebras, their
$(i,j)$ entries will have the parity $\bi+\bj\mod 2$.
The algebra of
even matrices over a superalgebra $\Ac$ will be identified with the tensor product algebra
$\End\CC^{2n+1|2m}\ot\Ac$, so that a matrix $A=[a_{ij}]$ is regarded as the element
\ben
A=\sum_{i,j=1}^{1'}e_{ij}\ot a_{ij}(-1)^{\bi\tss\bj+\bj}\in \End\CC^{2n+1|2m}\ot\Ac.
\een
We will use the involutive matrix {\em super-transposition} $t$ defined by
$(A^t)_{ij}=a_{j'i'}(-1)^{\bi\bj+\bj}\tss\ta_i\ta_j$,
where for $i=1,\dots,1'$ we set
\ben
\ta_i=\begin{cases} -1\qquad\text{if}\quad i>m+n+1\fand \bi=1,\\
\phantom{-}1\qquad\text{otherwise.}
\end{cases}
\een
We will also regard $t$ as the linear map
\beql{suptra}
t:\End\CC^{2n+1|2m}\to \End\CC^{2n+1|2m}, \qquad
e_{ij}\mapsto e_{j'i'}(-1)^{\bi\bj+\bi}\tss\ta_i\ta_j.
\eeq
In the case of multiple tensor products of the endomorphism algebras,
we will indicate by $t_a$ the map \eqref{suptra}
acting on the $a$-th copy of $\End\CC^{2n+1|2m}$.

A standard basis of the general linear Lie superalgebra $\gl_{2n+1|2m}$ is formed by elements $E_{ij}$
of the parity $\bi+\bj\mod 2$ for $1\leqslant i,j\leqslant 1'$ with the commutation relations
\ben
[E_{ij},E_{kl}]
=\de_{kj}\ts E_{i\tss l}-\de_{i\tss l}\ts E_{kj}(-1)^{(\bi+\bj)(\bk+\bl)}.
\een
We will regard the orthosymplectic Lie superalgebra $\osp_{2n+1|2m}$
as the subalgebra
of $\gl_{2n+1|2m}$ spanned by the elements
\ben
F_{ij}=E_{ij}-E_{j'i'}(-1)^{\bi\tss\bj+\bi}\ts\ta_i\ta_j.
\een
Introduce the permutation operator $P$ by
\ben
P=\sum_{i,j=1}^{1'} e_{ij}\ot e_{ji}(-1)^{\bj}\in \End\CC^{2n+1|2m}\ot\End\CC^{2n+1|2m}
\een
and set
\ben
Q=P^{\tss t_1}=P^{\tss t_2}=\sum_{i,j=1}^{1'} e_{ij}\ot e_{i'j'}(-1)^{\bi\bj}\ts\ta_i\ta_j
\in \End\CC^{2n+1|2m}\ot\End\CC^{2n+1|2m}.
\een
The $R$-{\em matrix} associated with $\osp_{2n+1|2m}$ is the
rational function in $u$ given by
\ben
R(u)=1-\frac{P}{u}+\frac{Q}{u-\ka},\qquad \ka=n-m-1/2.
\een
This is a super-version of the $R$-matrix
originally found in \cite{zz:rf}.
Following \cite{aacfr:rp}, we
define the {\it extended Yangian\/}
$\X(\osp_{2n+1|2m})=\X(\osp^{\tss\se}_{2n+1|2m})$ corresponding to the parity sequence $\se$
as a $\ZZ_2$-graded algebra with generators
$t_{ij}^{(r)}$ of parity $\bi+\bj\mod 2$, where $1\leqslant i,j\leqslant 1'$ and $r=1,2,\dots$,
satisfying the following defining relations.
Introduce the formal series
\beql{tiju}
t_{ij}(u)=\de_{ij}+\sum_{r=1}^{\infty}t_{ij}^{(r)}\ts u^{-r}
\in\X(\osp_{2n+1|2m})[[u^{-1}]]
\eeq
and combine them into the matrix $T(u)=[t_{ij}(u)]$.
Consider the elements of the tensor product algebra
$\End\CC^{2n+1|2m}\ot\End\CC^{2n+1|2m}\ot \X(\osp_{2n+1|2m})[[u^{-1}]]$ given by
\ben
T_1(u)=\sum_{i,j=1}^{1'} e_{ij}\ot 1\ot t_{ij}(u)(-1)^{\bi\tss\bj+\bj}\fand
T_2(u)=\sum_{i,j=1}^{1'} 1\ot e_{ij}\ot t_{ij}(u)(-1)^{\bi\tss\bj+\bj}.
\een
The defining relations for the algebra $\X(\osp_{2n+1|2m})$ take
the form of the $RTT$-{\em relation}
\beql{RTT}
R(u-v)\ts T_1(u)\ts T_2(v)=T_2(v)\ts T_1(u)\ts R(u-v).
\eeq
The algebra $\X(\osp^{\tss\se}_{2n+1|2m})$ associated with an arbitrary parity sequence $\se$
is isomorphic to the Yangian $\X(\osp_{2n+1|2m})$ associated with the standard sequence $\se^{\st}$;
an isomorphism is given by the map
$
t_{ij}(u)\mapsto t_{\si(i),\si(j)}(u)
$
for a suitable permutation $\si$ of the set $\{1,\dots,1'\}$.

The product $T(u-\ka)\ts T^{\tss t}(u)$ is a scalar matrix with
\beql{ttra}
T(u-\ka)\ts T^{\tss t}(u)=c(u)\tss 1,
\eeq
where $c(u)$ is a series in $u^{-1}$. Moreover, all the coefficients of the series $c(u)$ belong to
the center $\ZX(\osp_{2n+1|2m})$ of $\X(\osp_{2n+1|2m})$ and generate the center.

The {\em Yangian} $\Y(\osp_{2n+1|2m})$
is defined as the subalgebra of
$\X(\osp_{2n+1|2m})$ which
consists of the elements stable under
the automorphisms
\beql{muf}
t_{ij}(u)\mapsto f(u)\ts t_{ij}(u)
\eeq
for all series
$f(u)\in 1+u^{-1}\CC[[u^{-1}]]$.
We have the tensor product decomposition
\beql{tensordecom}
\X(\osp_{2n+1|2m})=\ZX(\osp_{2n+1|2m})\ot \Y(\osp_{2n+1|2m}).
\eeq
The Yangian $\Y(\osp_{2n+1|2m})$ can also be regarded as the quotient
of $\X(\osp_{2n+1|2m})$
by the relation $c(u)=1$.

Explicitly, the defining relations \eqref{RTT} can be written
with the use of super-commutator in terms of the series \eqref{tiju} as follows:
\begin{align}
\big[\tss t_{ij}(u),t_{kl}(v)\big]&=\frac{1}{u-v}
\big(t_{kj}(u)\ts t_{il}(v)-t_{kj}(v)\ts t_{il}(u)\big)
(-1)^{\bi\tss\bj+\bi\tss\bk+\bj\tss\bk}
\non\\
{}&-\frac{1}{u-v-\ka}
\Big(\de_{k i\pr}\sum_{p=1}^{1'}\ts t_{pj}(u)\ts t_{p'l}(v)
(-1)^{\bi+\bi\tss\bj+\bj\tss\bp}\ts\ta_i\ta_p
\label{defrel}\\
&\qquad\qquad\quad
{}-\de_{l j\pr}\sum_{p=1}^{1'}\ts t_{k\tss p'}(v)\ts t_{ip}(u)
(-1)^{\bi\tss\bk+\bj\tss\bk+\bi\tss\bp}\ts\ta_{j'}\ta_{p'}\Big).
\non
\end{align}

Each of the mappings
\beql{shift}
t_{ij}(u)\mapsto t_{ij}(u-a),\quad a\in \CC,
\eeq
and
\beql{transpo}
t_{ij}(u)\mapsto t_{j'i'}(-u)(-1)^{\bi\bj+\bj}\ts\ta_i\ta_j
\eeq
defines an automorphism of $\X(\osp_{2n+1|2m})$.
The universal enveloping algebra $\U(\osp_{2n+1|2m})$ can be regarded as a subalgebra of
$\X(\osp_{2n+1|2m})$ via the embedding
\beql{embosp}
F_{ij}\mapsto \frac12\big(t_{ij}^{(1)}-t_{j'i'}^{(1)}(-1)^{\bi\bj+\bj}\ts\ta_i\ta_j\big)(-1)^{\bi}.
\eeq
This fact relies on the Poincar\'e--Birkhoff--Witt theorem for the orthosymplectic Yangian
which was pointed out in \cite{aacfr:rp} and a detailed proof is given in \cite{gk:yo}; cf.
\cite[Sec.~3]{amr:rp}. It states that the algebra $\X(\osp_{2n+1|2m})$ is generated by
the coefficients of the series $c(u)$ and $t_{ij}(u)$ with the conditions
\ben
\bal
i+j&\leqslant 2n+2m+2\qquad \text{for}\quad i=1,\dots,m,m',\dots,1'\fand\\
i+j&< 2n+2m+2\qquad \text{for}\quad i=m+1,\dots,(m+1)'.
\eal
\een
Moreover, given any total ordering
on the set of the generators, the ordered monomials with the powers of odd generators
not exceeding $1$, form a basis of the algebra.

For a given parity sequence $\se=0\tss\se'$ which begins with $0$ consider
the extended Yangian $\X(\osp^{\tss\se'}_{2n-1|2m})$ and
for the parity sequence $\se=1\tss\se'$ beginning with $1$ consider
the extended Yangian $\X(\osp^{\tss\se'}_{2n+1|2m-2})$. In both cases, let
the indices
of the generators $t_{ij}^{(r)}$ of these algebras range over the sets
$2\leqslant i,j\leqslant 2\pr$ and $r=1,2,\dots$. The following
embedding properties were proved in \cite[Thm~3.1]{m:dt} for a standard
parity sequence, and the proof expends to arbitrary sequences $\se$ without any
significant changes. The mapping
\beql{embedgen}
t_{ij}(u)\mapsto t_{ij}(u)-t_{i1}(u)\ts t_{11}(u)^{-1}\tss t_{1j}(u),\qquad 2\leqslant i,j\leqslant 2\pr,
\eeq
defines the injective
homomorphisms
\beql{emb}
\X(\osp^{\tss\se'}_{2n-1|2m})\hra \X(\osp^{\tss\se}_{2n+1|2m})
\Fand
\X(\osp^{\tss\se'}_{2n+1|2m-2})\hra \X(\osp^{\tss\se}_{2n+1|2m})
\eeq
in the cases
$\se=0\tss\se'$ and $\se=1\tss\se'$, respectively.

The extended Yangian $\X(\osp_{2n+1|2m})$ is a Hopf algebra with the coproduct
defined by
\beql{Delta}
\De: t_{ij}(u)\mapsto \sum_{k=1}^{1'} t_{ik}(u)\ot t_{kj}(u).
\eeq
The image of the series $c(u)$ is found by the relation
$\De:c(u)\mapsto c(u)\ot c(u)$ and so the Yangian
$\Y(\osp_{2n+1|2m})$ inherits the Hopf algebra structure from $\X(\osp_{2n+1|2m})$.

\section{Highest weight representations}
\label{sec:hw}

We will keep a parity sequence $\se$ fixed.
A representation $V$ of the algebra $\X(\osp^{\tss\se}_{2n+1|2m})$
is called a {\em highest weight representation}
if there exists a nonzero vector
$\xi\in V$ such that $V$ is generated by $\xi$,
\begin{alignat}{2}
t_{ij}(u)\ts\xi&=0 \qquad &&\text{for}
\quad 1\leqslant i<j\leqslant 1', \qquad \text{and}\non\\
t_{ii}(u)\ts\xi&=\la_i(u)\ts\xi \qquad &&\text{for}
\quad i=1,\dots,1',
\label{trianb}
\end{alignat}
for some formal series
\beql{laiu}
\la_i(u)\in 1+u^{-1}\CC[[u^{-1}]].
\eeq
The vector $\xi$ is called the {\em highest vector}
of $V$.

\bpr\label{prop:nontrvm}
The series $\la_i(u)$ associated with a highest weight representation $V$
satisfy
the consistency conditions
\begin{multline}\label{nontrvm}
\la_i(u)\tss \la_{i\pr}\big(u-n+m+1/2+(-1)^{\bar 1}+\dots+(-1)^{\bi}\big)\\
{}=\la_{i+1}(u)\tss \la_{(i+1)'}\big(u-n+m+1/2+(-1)^{\bar 1}+\dots+(-1)^{\bi}\big)
\end{multline}
for $i=1,\dots,m+n$.
Moreover, the coefficients of the series $c(u)$ act in the representation
$V$ as the multiplications by scalars
determined by
$
c(u)\mapsto \la_1(u)\tss \la_{1'}(u-n+m+1/2).
$
\epr

\bpf
Relation \eqref{nontrvm} for $i=1$ follows by applying \eqref{defrel} to calculate
$t_{12}(u)\ts t_{1'2'}(v)\ts \xi$ and then setting
$v=u-\ka+(-1)^{\bar 1}$; cf. \cite[Prop.~4.4]{m:ry}. The formula
for the remaining values of $i$ is then implied by the embedding properties \eqref{emb}.
The eigenvalues of the coefficients of the series
$c(u)$ are found by taking the $(1',1')$ entry in
the matrix relation \eqref{ttra}.
\epf

As Proposition~\ref{prop:nontrvm} shows, the series $\la_i(u)$ in \eqref{trianb}
with $i>m+n+1$ are uniquely
determined by the first $m+n+1$ series. We will call the corresponding tuple
$\la(u)=(\la_{1}(u),\dots,\la_{m+n+1}(u))$
the {\em highest weight\/} of $V$.

Given an arbitrary tuple $\la(u)=(\la_{1}(u),\dots,\la_{m+n+1}(u))$
of formal series of the form \eqref{laiu}, define
the {\em Verma module} $M(\la(u))$ as the quotient of the algebra $\X(\osp^{\tss\se}_{2n+1|2m})$ by
the left ideal generated by all coefficients of the series $t_{ij}(u)$
with $1\leqslant i<j\leqslant 1'$ and $t_{ii}(u)-\la_i(u)$ for
$i=1,\dots,m+n+1$. The Poincar\'e--Birkhoff--Witt theorem for the algebra $\X(\osp^{\tss\se}_{2n+1|2m})$
implies that the Verma module $M(\la(u))$
is nonzero, and we denote by $L(\la(u))$ its irreducible quotient.
It is clear that the isomorphism
class of $L(\la(u))$ is determined by $\la(u)$.

\bpr\label{prop:fdhw}
Every finite-dimensional irreducible representation of the extended
Yangian $\X(\osp^{\tss\se}_{2n+1|2m})$
is isomorphic to the highest weight representation $L(\la(u))$ for a certain highest weight
$\la(u)=(\la_{1}(u),\dots,\la_{m+n+1}(u))$.
\epr

\bpf
The proof follows by a standard argument;
cf. \cite[Thm~5.1]{amr:rp}.
\epf

\subsection{Proof of Theorem~\ref{thm:necess}}
\label{subsec:thmnec}

As we outlined in the Introduction, we will use the embeddings \eqref{embedya}.
Recall from \cite{n:qb} that the Yangian $\Y(\gl_{n|m})$
(corresponding to the standard parity sequence $\se^{\st}$) is defined as the $\ZZ_2$-graded algebra
with generators $\bar t_{ij}^{\ts(r)}$
of parity $\bi+\bj\mod 2$, where $1\leqslant i,j\leqslant m+n$ and $r=1,2,\dots$,
while
\ben
\bi=\begin{cases} 1\qquad\text{for}\quad i=1,\dots,m,\\
0\qquad\text{for}\quad i=m+1,\dots,m+n.
\end{cases}
\een
The defining relations can be written in terms of the generating series
\ben
\bar t_{ij}(u)=\de_{ij}+\sum_{r=1}^{\infty}\bar t_{ij}^{\ts(r)}\ts u^{-r}
\in\Y(\gl_{n|m})[[u^{-1}]]
\een
and they have the form
\ben
\big[\tss \bar t_{ij}(u),\bar t_{kl}(v)\big]=\frac{1}{u-v}
\big(\bar t_{kj}(u)\ts \bar t_{il}(v)-\bar t_{kj}(v)\ts \bar t_{il}(u)\big)
(-1)^{\bi\tss\bj+\bi\tss\bk+\bj\tss\bk}.
\een
The first embedding in \eqref{embedya} is given by
\ben
\bar t_{ij}(u)\mapsto t_{ij}(u),\qquad 1\leqslant i,j\leqslant m+n.
\een
Since the representation $L(\la(u))$ is finite-dimensional, then so is
the $\Y(\gl_{n|m})$-module $\Y(\gl_{n|m})\xi$ defined via this embedding.
This is a highest weight representation with the highest weight
$(\la_1(u),\dots,\la_{m+n}(u))$. Hence, conditions \eqref{condgl} and \eqref{condo}
(the latter excluding $j=m+n$) must hold due to the classification theorem of \cite{zh:sy}
(see also its review in \cite{m:or} in the context of odd reflections).

The second embedding in \eqref{embedya} was given in \cite[Cor.~3.2]{m:dt}.
The same argument as with the first embedding shows that
since the cyclic span $\X(\oa_{2n+1})\xi$ is a finite-dimensional highest weight
representation with the highest weight
$(\la_{m+1}(u),\dots,\la_{m+n+1}(u))$, conditions \eqref{condo}
must hold due to the classification theorem of
\cite{amr:rp} in the case of the Yangian $\X(\oa_{2n+1})$.

Conditions \eqref{condgl} and \eqref{condo} imply that
by twisting
the representation $L(\la(u))$ by a suitable automorphism \eqref{muf},
if necessary, we may assume that all components $\la_i(u)$ of
the highest weight $\la(u)$ are polynomials in $u^{-1}$; that is,
for some $p\in\ZZ_+$ we have
\beql{hwdec}
\la_i(u)=(1+\la^{(1)}_i u^{-1})\dots (1+\la^{(p)}_i u^{-1}),\qquad i=1,\dots,m+n+1,
\eeq
with $\la^{(r)}_i\in\CC$ for $r=1,\dots,p$.

Consider the parity sequence $\se$ obtained from $\se^{\st}$ by replacing
the subsequence $\se_m\se_{m+1}=10$ with $01$. To calculate
the highest weight of the module $L(\la(u))$ associated with $\se$,
apply the corresponding odd reflection by using \cite[Thm~4.4]{m:or}.
We will assume that the parameters $\la^{(r)}_m$ and $\la^{(r)}_{m+1}$
are numbered in such a way that
\beql{laeq}
\la_m^{(r)}=\la_{m+1}^{(r)}\qquad \text{for all}\quad r=k+1,\dots,p
\eeq
for certain $k\in\{0,1,\dots,p\}$, while $\la_m^{(r)}\ne\la_{m+1}^{(s)}$
for all $1\leqslant r,s\leqslant k$. Then set
\ben
\bal
\la^{[1]}_m(u)&=(1+(\la_{m}^{(1)}+1)\tss u^{-1})\dots (1+(\la_{m}^{(k)}+1)\tss u^{-1})
(1+\la_{m}^{(k+1)}u^{-1})\dots (1+\la_{m}^{(p)}u^{-1}),\\[0.4em]
\la^{[1]}_{m+1}(u)&=(1+(\la_{m+1}^{(1)}+1)\tss u^{-1})\dots (1+(\la_{m+1}^{(k)}+1)\tss u^{-1})
(1+\la_{m}^{(k+1)}u^{-1})\dots (1+\la_{m}^{(p)}u^{-1}).
\eal
\een
Note that these series coincide with those given by \eqref{lammpo}.

Although the arguments of \cite{m:or} deal with
representations of $\Y(\gl_{n|m})$, the odd reflections apply
to representations of the Yangian $\X(\osp_{2n+1|2m})$ as well, by taking into account the
first embedding in \eqref{embedya}. The final step in the proof of Theorem~4.4 therein
should just be adjusted to use the re-labelling automorphism of $\X(\osp_{2n+1|2m})$
acting on the generators by
\ben
t_{ij}(u)\mapsto t_{\si(i),\si(j)}(u),
\een
where $\si\in\Sym_{2n+2m+1}$ is the product of the transpositions $(m,m+1)$ and $((m+1)',m')$.

As the next step, apply the odd reflection to the pair $(\la^{[1]}_m(u),\la_{m+2}(u))$
and then continue as described in \eqref{oddlami} to conclude that the highest weight
of the module $L(\la(u))$ associated with the parity sequence $\se=1\dots 10\dots 01$
has the form
\ben
\big(\la_1(u),\dots,\la_{m-1}(u),\la^{[1]}_{m+1}(u),\dots,\la^{[1]}_{m+n}(u),
\la^{[n]}_{m}(u),\la_{m+n+1}(u)\big).
\een

Finally, we use a composition of embeddings \eqref{embedgen} as in
\cite[Cor.~3.2]{m:dt} to get the embedding given in \eqref{embospyang}.
The cyclic span $\X(\osp_{1|2})\xi$ is a finite-dimensional highest weight
representation with the highest weight
$(\la^{[n]}_{m}(u),\la_{m+n+1}(u))$. Thus conditions \eqref{condosp}
follow from the Main Theorem of \cite{m:ry}.

\section{Representations with linear highest weights}
\label{sec:rl}

Here we will be concerned with the representations $L(\la(u))$ of the Yangian
$\X(\osp_{2n+1|2m})$ (associated with the standard parity sequence $\se^{\st}$), such that
all components of the highest weight $\la(u)$ are linear in $u^{-1}$
as given in \eqref{linhw}. Observe that the twist of $L(\la(u))$ by the composition
of the automorphism \eqref{muf} for $f(u)=u/(u+\la_{m+n+1})$
and the automorphism \eqref{shift} for $a=-\la_{m+n+1}$ yields a representation
with the linear highest weight, where the components are changed by the rule
$\la_i\mapsto \la_i-\la_{m+n+1}$. Therefore, we will not restrict generality
by assuming $\la_{m+n+1}=0$ so that
\beql{linhwlaze}
\la(u)=(1+\la_1\ts u^{-1},\dots,1+\la_{m+n}\ts u^{-1},1).
\eeq

The cyclic span $\U(\osp_{2n+1|2m})\xi$
defined via the embedding \eqref{embosp} is a highest weight representation
of $\osp_{2n+1|2m}$ with $F_{ij}\tss\xi=0$ for all $1\leqslant i<j\leqslant 1'$.
It follows from \eqref{ttra} that the generator $t^{(1)}_{m+n+1,m+n+1}$
belongs to the center of $\X(\osp_{2n+1|2m})$ and so it acts as multiplication
by $0$ in $L(\la(u))$. Hence, $F_{ii}\tss\xi=(-1)^{\bi}\tss\la_i\tss\xi$ for $i=1,\dots,m+n$
by \eqref{embosp},
and so
the highest weight
of the $\osp_{2n+1|2m}$-module $\overline L=\U(\osp_{2n+1|2m})\xi$ is given by
\beql{osphwls}
(-\la_1,\dots,-\la_m,\la_{m+1},\dots,\la_{m+n}).
\eeq

For complex numbers $a$ and $b$ we will write $a\rar b$ or $b\lar a$ to mean
that $a-b\in\ZZ_+$.

\bth\label{thm:linear}
The representation $L(\la(u))$ of $\X(\osp_{2n+1|2m})$ with the highest weight \eqref{linhwlaze}
is finite-dimensional if and only if either
\beql{halfint}
\la_1\lar\la_2\lar\cdots\lar\la_m\lar -n\Fand \la_{m+1}\rar\cdots\rar\la_{m+n}\rar 1/2,
\eeq
or
\beql{ydiag}
\la_1\lar\cdots\lar\la_m=-l\fand \la_{m+1}\rar\cdots\rar\la_{m+l}\rar 0=\la_{m+l+1}=\cdots=\la_{m+n}
\eeq
for some $l\in\ZZ_+$; if $l\geqslant n$ then the second part of condition \eqref{ydiag}
is understood as
\ben
\la_{m+1}\rar\cdots\rar\la_{m+n}\rar 0.
\een
\eth

\bpf
Suppose first that $\dim L(\la(u))<\infty$. Since the cyclic span
$\overline L$ is finite-dimensional, the $\osp_{2n+1|2m}$-highest weight
\eqref{osphwls} must satisfy the finite-dimensionality conditions for irreducible
highest weight representations of $\osp_{2n+1|2m}$ as obtained by Kac~\cite{k:rc};
see also \cite[Sec.~2.1]{cw:dr}. This implies the desired conditions on the
components of the highest weight.

Conversely, suppose that the conditions on the
components of the highest weight $\la(u)$ are satisfied.
Equip the set of $\osp_{2n+1|2m}$-weights of
$\overline L$
with the standard partial
order determined by
the choice of the positive root system corresponding to the upper-triangular
subalgebra of $\osp_{2n+1|2m}$.
The next lemma will be our main tool for proving that certain vectors in
$L(\la(u))$ are zero.

\ble\label{lem:uppze}
Suppose that $\eta\in L(\la(u))$ is
an $\osp_{2n+1|2m}$-weight vector whose weight is lower than the weight of $\xi$
and let $S$ be a subset of $\{1,\dots,m+n\}$.
Then the conditions
\beql{coanni}
t_{i, i+1}(u)\tss\eta=0\quad\text{for}\quad i\in S\fand
t_{(i+1)', i'}(u) \tss\eta=0\quad\text{for}\quad i\in \{1,\dots,m+n\}\setminus S
\eeq
imply that $\eta=0$.
\ele

\bpf
Suppose first that $S=\{1,\dots,m+n\}$. By taking commutators of both sides of
the relations $t_{i, i+1}(u)\tss\eta=0$
with suitable elements $t^{(1)}_{j, j+1}$ with $j=1,\dots,m+n$ we come
to the conditions
\beql{condva}
t_{ij}(u)\tss \eta=0\qquad\text{for all}\quad 1\leqslant i<j\leqslant m+n+1.
\eeq
Since the action of the coefficients
of the series $t_{jj}(u)$ on the vector $\eta$ does not change its $\osp_{2n+1|2m}$-weight,
it follows from
the Poincar\'e--Birkhoff--Witt theorem that a nonzero vector $\eta$ satisfying \eqref{condva}
would generate a submodule of $L(\la(u))$ which does not contain $\xi$.
However, this is impossible due to the irreducibility of $L(\la(u))$, thus proving
that $\eta=0$.

The proof for an arbitrary subset $S\in\{1,\dots,m+n\}$ follows from the observation
that the condition $t_{i, i+1}(u)\tss\eta=0$ for any given $i\in\{1,\dots,m+n\}$ can be replaced with
$t_{(i+1)', i'}(u)\tss\eta=0$ for the conclusion to remain valid.
This is clear from relation \eqref{ttra}
written in the form
\ben
T^{\tss t}(u)=c(u)\tss T(u-\ka)^{-1}.
\een
By using the expansion
\ben
T(u)=1+T^{(1)}u^{-1}+T^{(2)}u^{-2}+\dots
\een
and setting the degree of $t_{ij}^{(r)}$ equal to $r$, we derive that
\ben
t_{(i+1)', i'}^{(r)}=\pm\tss t_{i, i+1}^{(r)}\quad+\quad\text{a linear combination of monomials of degree}\ <r.
\een
The claim is now verified by an easy induction on $r$.
\epf

We proceed by considering a particular representation which will play the role
of a fundamental Yangian module.

\ble\label{lem:fundosc}
The representation $L^{\circ}=L(\la(u))$ of $\X(\osp_{2n+1|2m})$ with the highest weight
\eqref{linhwlaze} with
\ben
\la_1=\dots=\la_m=-n\Fand \la_{m+1}=\dots=\la_{m+n}=1/2
\een
is finite-dimensional.
\ele

\bpf
We will be using the eigenvalues $\la_{i'}(u)$
of the diagonal operators $t_{i'i'}(u)$ on $\xi$ for
the values $i=1,\dots,m+n$. They are found in Proposition~\ref{prop:nontrvm}
and given by
\ben
\la_{i'}(u)=1+n\tss u^{-1}\qquad\text{for}\quad i=1,\dots,m
\een
and
\ben
\la_{i'}(u)=1-\frac12\tss u^{-1}\qquad\text{for}\quad i=m+1,\dots,m+n.
\een

The highest weight \eqref{osphwls} associated with the given
components $\la_i$
corresponds to a finite-dimensional representation of the Lie superalgebra
$\osp_{2n+1|2m}$. Therefore, it will be sufficient to demonstrate
that all elements $t_{ij}^{(r)}$ of the extended Yangian with $r\geqslant 2$
act on $\overline L$ as the zero operators. Furthermore, it is clear from
the defining relations \eqref{defrel} that it will be enough to
verify that these elements act as the zero operators on the highest vector $\xi$.

Observe that
\beql{toopo}
t_{j+1,\tss j}(v)\tss\xi=0\qquad\text{for}\quad j=1,\dots,m-1,m+1,\dots,m+n-1.
\eeq
Indeed, by \eqref{defrel}, up to a sign factor depending on $j$, the expression
$t_{j,\tss j+1}(u)\ts t_{j+1,\tss j}(v)\tss\xi$ equals
\ben
\frac{1}{u-v}
\big(t_{j+1,\tss j+1}(u)\ts t_{jj}(v)-t_{j+1,\tss j+1}(v)\ts t_{jj}(u)\big)\tss\xi
\een
which vanishes because $\la_j(u)=\la_{j+1}(u)$. Furthermore,
$t_{j+1,\tss j}(v)\tss\xi$ is annihilated by the remaining simple root series
$t_{i,i+1}(u)$ with $i\ne j$. This follows by applying \eqref{defrel}
to $[t_{i,i+1}(u),t_{j+1,\tss j}(v)]$ for $i<j$, and to $[t_{j+1,\tss j}(v),t_{i,i+1}(u)]$
for $i>j$. Thus \eqref{toopo} follows
from Lemma~\ref{lem:uppze}.

By taking commutators of the left hand side of \eqref{toopo} with
suitable elements $t^{(1)}_{i+1, i}$ we come to more general relations
\beql{toopoge}
t_{j\ts i}(v)\tss\xi=0\qquad\text{for}\quad 1\leqslant i<j\leqslant m
\fand m+1\leqslant i<j\leqslant m+n.
\eeq
Since the transposition automorphism \eqref{transpo} preserves the module $L(\la(u))$,
relations \eqref{toopoge} extend to the values
$1\leqslant j'<i'\leqslant m$ and $m+1\leqslant j'<i'\leqslant m+n$.

As the next step, we verify that
\beql{tmpom}
\big(t_{m+1,m}(v)-v^{-1}\tss t^{(1)}_{m+1,m}\big)\ts\xi=0.
\eeq
It is clear from \eqref{defrel} that the extression
on the left hand side is annihilated by the series $t_{i, i+1}(u)$
for $i=1,\dots,m-1,m+1,\dots,m+n$. Furthermore, by \eqref{defrel}
the expression
\ben
t_{m, m+1}(u)\big(t_{m+1,m}(v)-v^{-1}\tss t^{(1)}_{m+1,m}\big)\ts\xi
\een
equals
\ben
\frac{1}{u-v}
\big(t_{m+1, m+1}(u)\ts t_{mm}(v)-t_{m+1, m+1}(v)\ts t_{mm}(u)\big)\tss\xi
+\frac{1}{v}\ts\big(t_{m+1, m+1}(u)-t_{mm}(u)\big)\tss\xi.
\een
The coefficient of $\xi$ simplifies to
\ben
\frac{1}{u-v}\Big(\frac{(u+1/2)(v-n)}{uv}-\frac{(v+1/2)(u-n)}{uv}\Big)
+\frac{1/2+n}{uv}=0
\een
and \eqref{tmpom} follows from Lemma~\ref{lem:uppze}.

Take commutators of the left hand side of \eqref{tmpom} with
suitable elements $t^{(1)}_{i+1\ts i}$, where $i$ satisfies the conditions
of \eqref{toopo}, to derive that
\beql{tmpomge}
\big(t_{m+j, i}(v)-v^{-1}\tss t^{(1)}_{m+j, i}\big)\ts\xi=0\qquad\text{for}\quad i=1,\dots,m
\fand j=1,\dots,n.
\eeq
The same calculation verifies that
\ben
\big(t_{m+n+1,m+n}(v)-v^{-1}\tss t^{(1)}_{m+n+1,m+n}\big)\ts\xi=0
\een
and hence
\beql{tmomn}
\big(t_{m+n+1,m+j}(v)-v^{-1}\tss t^{(1)}_{m+n+1,m+j}\big)\ts\xi=0
\qquad\text{for}\quad j=1,\dots,n.
\eeq
The application of the automorphism \eqref{transpo} yields the respective counterparts
of \eqref{tmpomge} and \eqref{tmomn} for the transposed series.

Furthermore, to derive that
\ben
\big(t_{m+n+1,m}(v)-v^{-1}\tss t^{(1)}_{m+n+1,m}\big)\ts\xi=0
\een
we note first that the expression on the left hand side is annihilated
by the series $t_{i, i+1}(u)$
for $i=1,\dots,m$ as verified by the same calculation as for \eqref{tmpom}.
For $i=m+1,\dots,m+n$ the expression
\ben
t_{i, i+1}(u)\ts \big(t_{m+n+1,m}(v)-v^{-1}\tss t^{(1)}_{m+n+1,m}\big)\ts\xi
\een
equals
\ben
\frac{1}{u-v}
\big(t_{m+n+1, i+1}(u)\ts t_{im}(v)-t_{m+n+1, i+1}(v)\ts t_{im}(u)\big)\tss\xi
+\frac{1}{v}\ts\de_{i,m+n}\tss t_{im}(u)\tss\xi.
\een
The previously verified relations \eqref{tmpomge} and \eqref{toopoge} imply that
this expression is zero. Hence, we may conclude that
\ben
\big(t_{m+n+1,i}(v)-v^{-1}\tss t^{(1)}_{m+n+1,i}\big)\ts\xi=0
\qquad\text{for}\quad i=1,\dots,m+n
\een
together with the transposed relations implies by the application
of the automorphism \eqref{transpo}.

Finally, we will show by a reverse induction that for each $j=1,\dots,m+n$ we have
\beql{reain}
\big(t_{j'i}(v)-v^{-1}\tss t^{(1)}_{j'i}\big)\ts\xi=0 \fand
\big(t_{i'j}(v)-v^{-1}\tss t^{(1)}_{i'j}\big)\ts\xi=0
\qquad\text{for}\quad i=1,\dots,j.
\eeq
Taking first $i=j$ note the relation
\ben
\big(t_{j'j}(v)-v^{-1}\tss t^{(1)}_{j'j}\big)\ts\xi
=\big(t_{j',\tss j+1}(v)-v^{-1}\tss t^{(1)}_{j',\tss j+1}\big)\ts t^{(1)}_{j+1,\tss j}\xi
\een
which holds by the induction hypothesis. Hence, up to a sign factor depending on $j$,
the application of $t_{(j+1)',\tss j'}(u)$ to this expression yields
\beql{appjj}
\frac{1}{u-v}
\big(t_{j'j'}(u)\ts t_{(j+1)',\tss j+1}(v)-t_{j'j'}(v)
\ts t_{(j+1)',\tss j+1}(u)\big)\ts t^{(1)}_{j+1,\tss j}\xi
-\frac{1}{v}\ts t_{(j+1)',\tss j+1}(u)\ts t^{(1)}_{j+1,\tss j}\xi.
\eeq
By applying \eqref{defrel} and the induction hypothesis, we find that
\ben
t_{(j+1)',\tss j+1}(u)\ts t^{(1)}_{j+1,\tss j}\xi
=u^{-1}\tss t^{(1)}_{(j+1)',\tss j+1}\ts t^{(1)}_{j+1,\tss j}\xi.
\een
Furthermore, $t_{j'j'}(u)$ commutes with $t^{(1)}_{(j+1)',\tss j+1}$, whereas
\ben
t_{j'j'}(u)\ts t^{(1)}_{j+1,\tss j}\xi=\la_{j'}(u)\ts t^{(1)}_{j+1,\tss j}\xi
+\big[t_{j'j'}(u),\ts t^{(1)}_{j+1,\tss j}\big]\ts\xi=(1+\al\tss u^{-1})\ts t^{(1)}_{j+1,\tss j}\xi
\een
for a constant $\al$ depending on $j$. This implies that
expression \eqref{appjj} vanishes.

The application of the remaining series $t_{k,k+1}(u)$ with $k=1,\dots,j-1,j+1,\dots,m+n$
to the expression in \eqref{reain} with $i=j$ yields zero as well, which is seen by a
simpler calculation. Invoking Lemma~\ref{lem:uppze} again, we conclude that
the expression in \eqref{reain} with $i=j$ is zero.

The verification of the remaining relations in \eqref{reain}
is done by quite similar calculations and a subsequent use of the transposition
automorphism \eqref{transpo}.
\epf

We will point out the following super-analogue of the well-known property
of Yangian modules, which is immediate from the coproduct formula \eqref{Delta}.

\ble\label{lem:cophw}
Let $L(\mu(u))$ and $L(\nu(u))$ be the irreducible highest weight $\X(\osp_{2n+1|2m})$-modules
with the highest weights
\ben
\mu(u)=\big(\mu_1(u),\dots,\mu_{m+n+1}(u)\big)\Fand
\nu(u)=\big(\nu_1(u),\dots,\nu_{m+n+1}(u)\big)
\een
and the respective highest vectors $\eta$ and $\ze$.
Then the cyclic span $\X(\osp_{2n+1|2m})(\eta\ot\ze)$ is a highest weight module
with the highest weight
\beql{mhw}
\big(\mu_1(u)\nu_1(u),\dots,\mu_{m+n+1}(u)\nu_{m+n+1}(u)\big).
\eeq
Hence, if the modules $L(\mu(u))$ and $L(\nu(u))$ are finite-dimensional,
then so is the irreducible highest weight module
with the highest weight \eqref{mhw}.
\qed
\ele

We will need a family of fundamental modules over $\X(\osp_{2n+1|2m})$
which are analogous to polynomial representations of the Yangian $\Y(\gl_{n|m})$
as reviewed in \cite[Appendix]{m:rs}.

The {\em vector representation} of $\X(\osp^{\tss\se}_{2n+1|2m})$
(with an arbitrary parity sequence $\se$) on $\CC^{2n+1|2m}$ is defined by
\beql{vectre}
t_{ij}(u)\mapsto \de_{ij}+u^{-1}\tss
e_{ij}(-1)^{\bi}-(u+\ka)^{-1}\tss e_{j'i'}(-1)^{\bi\bj}\ts\ta_i\ta_j.
\eeq
The homomorphism property
follows from the $RTT$-relation \eqref{RTT} and
the Yang--Baxter equation satisfied by $R(u)$; cf. \cite{aacfr:rp}.

For the standard parity sequence $\se^{\st}$ and $d\in\{1,\dots,m\}$,
use the coproduct \eqref{Delta} to equip
the tensor product space $(\CC^{2n+1|2m})^{\ot d}$
with the action of $\X(\osp_{2n+1|2m})$ by setting
\ben
t_{ij}(u)\mapsto
\sum_{a_1,\dots,a_{d-1}=1}^{1'} t_{i,a_1}(u+d-1)\ot t_{a_1,a_2}(u+d-2)
\ot\dots\ot t_{a_{d-1},j}(u),
\een
where the generators act in the respective copies of the vector space
$(\CC^{2n+1|2m})^{\ot d}$ via the rule \eqref{vectre}.
Set
\beql{xid}
\xi_d=\sum_{\si\in\Sym_d} \sgn\si\cdot
e_{\si(1)}\ot\dots\ot e_{\si(d)} \in (\CC^{2n+1|2m})^{\ot d}.
\eeq
The calculations of \cite[Appendix]{m:rs} show that
the cyclic span $\X(\osp_{2n+1|2m})\ts \xi_d$ is a highest weight representation
of $\X(\osp_{2n+1|2m})$ with the highest weight $\la(u)$
whose components are found by
\ben
\la_i(u)=1-u^{-1}\quad\text{for}\quad i=1,\dots,d
\fand \la_i(u)=1\quad\text{for}\quad i=d+1,\dots,m+n+1.
\een
We will denote the irreducible quotient of this representation by $L^{\sharp\ts d}$.

Now consider the parity sequence $\se=0\dots01\dots1$
and recall that the algebra
$\X(\osp^{\tss\se}_{2n+1|2m})$ is isomorphic to the extended Yangian
$\X(\osp_{2n+1|2m})$ associated with the standard parity sequence.
For $d\in\{1,\dots,n\}$
equip
the tensor product space $(\CC^{2n+1|2m})^{\ot d}$
with the action of $\X(\osp^{\tss\se}_{2n+1|2m})$ by setting
\ben
t_{ij}(u)\mapsto
\sum_{a_1,\dots,a_{d-1}=1}^{1'} t_{i,a_1}(u-d+1)\ot t_{a_1,a_2}(u-d+2)
\ot\dots\ot t_{a_{d-1},j}(u),
\een
where the generators act in the respective copies of the vector space
$(\CC^{2n+1|2m})^{\ot d}$ via the rule \eqref{vectre}.
The vector $\xi_d$ defined by the same formula \eqref{xid}
now generates a highest weight representation
of $\X(\osp^{\tss\se}_{2n+1|2m})$ with the highest weight $\la(u)$
whose components are found by
\ben
\la_i(u)=1+u^{-1}\quad\text{for}\quad i=1,\dots,d
\fand \la_i(u)=1\quad\text{for}\quad i=d+1,\dots,m+n+1.
\een
We will denote the irreducible quotient of this representation by $L^{\flat\tss d}$.

It is clear from the description of the Yangian odd reflections in Sec.~\ref{subsec:thmnec},
that the linear highest weights of the form \eqref{linhwlaze}
are transformed in the same way as the highest weights
\eqref{osphwls} are transformed with respect to the corresponding
orthosymplectic Lie superalgebra odd reflections, as described
in \cite[Sec.~1.3]{cw:dr}. In particular, we get the following correspondence
between the linear highest weights of the modules associated with Young diagrams
implied by \cite[Sec.~2.4]{cw:dr}.

Recall that an $(m,n)$-{\em hook partition} $\Ga=(\Ga_1,\Ga_2,\dots)$ is a partition satisfying the
condition
$\Ga_{m+1}\leqslant n$. This means that the Young diagram $\Ga$ is contained
in the $(m,n)$-{\em hook} as depicted below. The figure also illustrates the partitions
$\mu=(\mu_1,\dots,\mu_m)$ and $\nu=(\nu_1,\dots,\nu_n)$
associated with $\Ga$. They are introduced by setting
\ben
\mu_i=\max\{\Ga_i-n,0\},\qquad i=1,\dots,m,
\een
and
\ben
\nu_j=\max\{\Ga\pr_j-m,0\},\qquad j=1,\dots,n,
\een
where $\Ga\pr$ denotes the conjugate partition so that $\Ga\pr_j$ is the length
of column $j$ in the diagram $\Ga$:

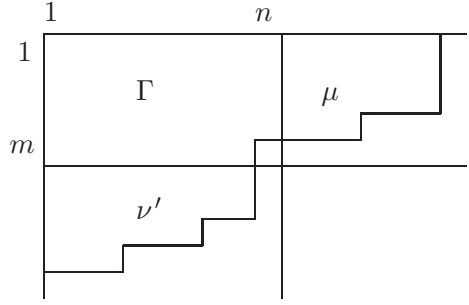
\begin{center}
\begin{picture}(150,90)
\thinlines

\put(0,0){\line(0,1){100}}
\put(90,0){\line(0,1){100}}
\put(0,50){\line(1,0){160}}
\put(0,100){\line(1,0){160}}

\put(0,10){\line(1,0){30}}
\put(30,10){\line(0,1){10}}
\put(30,20){\line(1,0){30}}
\put(60,20){\line(0,1){10}}
\put(60,30){\line(1,0){20}}
\put(80,30){\line(0,1){30}}
\put(80,60){\line(1,0){40}}
\put(120,60){\line(0,1){10}}
\put(120,70){\line(1,0){30}}
\put(150,70){\line(0,1){30}}

\put(0,105){\small$1$}
\put(80,105){\small$n$}
\put(-10,90){\small$1$}
\put(-13,55){\small$m$}

\put(35,75){\small$\Ga$}
\put(105,75){\small$\mu$}
\put(35,30){\small$\nu\pr$}

\end{picture}
\end{center}


\noindent
We will associate two $(m+n)$-tuples of integers with $\Ga$ by
\beql{ga}
\Ga^{\sharp}=(-\Ga_1,\dots,-\Ga_m,\nu_1,\dots,\nu_n)\Fand
\Ga^{\flat}=(\Ga'_1,\dots,\Ga'_n,-\mu_1,\dots,-\mu_m).
\eeq
According to \cite[Example~2.53]{cw:dr}, if the highest weight
\eqref{osphwls} for the standard parity sequence $\se^{\st}$ coincides with $\Ga^{\sharp}$ for
an $(m,n)$-hook partition $\Ga$, then the highest weight for
the parity sequence $\se=0\dots01\dots1$ obtained by a sequence of odd reflections,
coincides with $\Ga^{\flat}$.

We are now in a position to complete the proof of the theorem.
Suppose first that conditions \eqref{ydiag} hold.
They mean that
the highest weight \eqref{osphwls} coincides
with $\Ga^{\sharp}$ for certain $(m,n)$-hook partition $\Ga$.
Hence it will be sufficient to show that the highest module
$L(\la(u))$ associated with $\Ga$ in this way, is finite-dimensional.

Recall the irreducible highest weight representations of $L^{\sharp\tss d}$ and
$L^{\flat\tss d}$ of the extended Yangian constructed above and denote by
$L^{\sharp\tss d}_a$ and
$L^{\flat\tss d}_a$ their respective compositions with the shift
automorphism \eqref{shift}.

Given an $(m,n)$-hook partition $\Ga$, consider the tensor product module
\ben
L_{\Ga}=\bigotimes_{d=1}^m \Big(L^{\sharp\tss d}_{\Ga_d-1}\ot L^{\sharp\tss d}_{\Ga_d-2}\ot\dots\ot
L^{\sharp\tss d}_{\Ga_{d+1}}\Big),
\een
where $\Ga_{m+1}$ should be replaced by $0$. According to Lemma~\ref{lem:cophw}, the cyclic
$\X(\osp_{2n+1|2m})$-span of the tensor product of the highest weight vectors of
the modules $L^{\sharp\tss d}_{a}$ is a highest weight module with the highest weight
given by \eqref{linhwlaze} with
\ben
\la_i=-\Ga_i\quad\text{for}\quad i=1,\dots,m
\Fand \la_i=0\quad\text{for}\quad i=m+1,\dots,m+n.
\een
As we recalled above,
by applying a sequence of odd reflections, we find that
the highest weight of the irreducible quotient $\overline L_{\Ga}$ of this cyclic span
associated with the parity sequence $\se=0\dots01\dots1$ is found as
$\overline\Ga^{\tss\flat}$, where $\overline\Ga$ is the Young diagram
with $m$ rows $\Ga_1,\dots,\Ga_m$.

If the parameter $l$ in \eqref{ydiag} exceeds $n-1$, it should be understood as equal to $n$
in the argument below; in that case we set $\nu_{n+1}:=0$.
Consider
the tensor product module
\ben
\overline L_{\Ga}\ot\bigotimes_{d=1}^l \Big(L^{\flat\tss d}_{-\nu_d-m+1}
\ot L^{\flat\tss d}_{-\nu_d-m+2}\ot\dots\ot
L^{\flat\tss d}_{-\nu_{d+1}-m}\Big).
\een
By Lemma~\ref{lem:cophw}, the cyclic
$\X(\osp^{\tss\se}_{2n+1|2m})$-span of the tensor product of the highest weight vectors of
the tensor factors is a highest weight module with the highest weight
associated with $\Ga^{\flat}$. All modules
involved in the tensor products are finite-dimensional and so is the
irreducible quotient of the cyclic span. This completes the proof
of the sufficiency of conditions \eqref{ydiag}.

Now suppose that conditions \eqref{halfint} hold. The argument will be quite similar
to the above, taking the finite-dimensional representation $L^{\circ}$ constructed
in Lemma~\ref{lem:fundosc} as the starting point. Consider the tensor product module
\ben
L^{\circ}\ot\bigotimes_{d=1}^m \Big(L^{\sharp\tss d}_{-\la_d-1}\ot L^{\sharp\tss d}_{-\la_d-2}\ot\dots\ot
L^{\sharp\tss d}_{-\la_{d+1}}\Big),
\een
where we set $\la_{m+1}:=-n$. By Lemma~\ref{lem:cophw}, the cyclic
$\X(\osp_{2n+1|2m})$-span of the tensor product of the highest weight vectors of
the tensor factors is a highest weight module with the highest weight
given by
\ben
\big(1+\la_1\tss u^{-1},\dots,1+\la_m\tss u^{-1},1+1/2\tss u^{-1},\dots,1+1/2\tss u^{-1},1\big).
\een
By applying a sequence of odd reflections, we can regard
the irreducible quotient $\overline L^{\ts\circ}$
of the cyclic span as a $\X(\osp^{\tss\se}_{2n+1|2m})$-module with $\se=0\dots01\dots1$,
whose highest weight is found by
\ben
\big(1+(m+1/2)\tss u^{-1},\dots,1+(m+1/2)\tss u^{-1},
1+(\la_1+n)\tss u^{-1},\dots,1+(\la_m+n)\tss u^{-1},1\big).
\een
Finally, $L(\la(u))$, regarded as a $\X(\osp^{\tss\se}_{2n+1|2m})$-module, is isomorphic to
the irreducible quotient of the cyclic span of the tensor product of the highest weight vectors of
the tensor factors in
\ben
\overline L^{\ts\circ}\ot
\bigotimes_{d=1}^n \Big(L^{\flat\tss d}_{-\la_{m+d}-m+1}\ot L^{\flat\tss d}_{-\la_{m+d}-m+2}\ot\dots\ot
L^{\flat\tss d}_{-\la_{m+d+1}-m}\Big),
\een
where we set $\la_{m+n+1}:=1/2$. Thus, $L(\la(u))$ is finite-dimensional.
\epf

The following corollary confirms Conjecture~\ref{conj:fdim} in the case
of generic highest weights.

\bco\label{cor:genhw}
Suppose that the components of the highest weight $\la(u)$ are given
by \eqref{hwdec} with the condition that
for each $i=1,\dots,m+n+1$ none of the differences $\la^{(a)}_i-\la^{(b)}_i$
is an integer for $a\ne b$. Then the representation $L(\la(u))$ is finite-dimensional
if and only if $\la(u)$ satisfies the conditions of Theorem~\ref{thm:necess}.
\eco

\bpf
Due to Theorem~\ref{thm:necess}, we only need to establish the
sufficiency of the conditions. They imply that
because of the additional assumptions,
for each $i$
the parameters $\la^{(1)}_i,\dots,\la^{(p)}_i$ can be re-numbered in such a way
that each linear weight
\beql{splinhw}
\la^{(a)}(u)=\big(1+\la^{(a)}_1\tss u^{-1},\dots,1+\la^{(a)}_{m+n+1}\tss u^{-1}\big)
\eeq
satisfies the conditions of Theorem~\ref{thm:necess} for $a=1,\dots,p$. By Theorem~\ref{thm:linear},
each representation $L(\la^{(a)}(u))$ is finite-dimensional and so is
the cyclic span of the tensor products of the highest vectors in
\ben
L(\la^{(1)}(u))\ot\dots\ot L(\la^{(p)}(u)).
\een
By Lemma~\ref{lem:cophw}, the irreducible quotient of the cyclic span
is isomorphic to $L(\la(u))$, thus implying that it is finite-dimensional.
\epf

\section{Classification theorem for representations of $\X(\osp_{3|2m})$}
\label{sec:ct}

In this section we specialize to the case $n=1$ and prove Theorem~\ref{thm:none}.
We will show that the conditions of Theorem~\ref{thm:necess} imply that
the highest weight $\la(u)$ can be split into linear weights in a way similar
to the proof of Corollary~\ref{cor:genhw}.

As in the proof of Theorem~\ref{thm:necess},
by twisting
the representation $L(\la(u))$ by a suitable automorphism \eqref{muf},
we may assume that all components $\la_i(u)$ of
the highest weight $\la(u)$ are given by \eqref{hwdec}.
Furthermore, we may also
assume that the parameters $\la^{(r)}_m$ and $\la^{(r)}_{m+1}$
are numbered in such a way that relations \eqref{laeq} hold
for certain $k\in\{0,1,\dots,p\}$, while $\la_m^{(r)}\ne\la_{m+1}^{(s)}$
for all $1\leqslant r,s\leqslant k$. Then the application
of the odd reflection yields the subsequent formulas for
$\la^{[1]}_m(u)$ and $\la^{[1]}_{m+1}(u)$ given therein.

Suppose first that $\la_{m+2}^{(a)}$ is not equal to any of
$\la_m^{(k+1)},\dots,\la_m^{(p)}$ for each $a=1,\dots,p$.
The condition \eqref{condosp} reads
\ben
\frac{\la_{m+2}(u)}{\la^{[1]}_{m}(u)}=\frac{P_m(u+1)}{P_m(u)},
\een
and it
implies that the parameters $\la_{m+2}^{(a)}$ with $a=1,\dots,p$ can
be renumbered to satisfy the conditions
\beql{diffco}
\la_{m+2}^{(a)}-\la_m^{(a)}-1\in\ZZ_+\qquad\text{for}\qquad a=1,\dots,k
\eeq
and
\ben
\la_{m+2}^{(a)}-\la_m^{(a)}\in\ZZ_+\qquad\text{for}\qquad a=k+1,\dots,p.
\een
Since we excluded the equalities $\la_{m+2}^{(a)}=\la_m^{(a)}$,
we find that condition \eqref{diffco} extends to all values $a=1,\dots,p$.

Furthermore, condition \eqref{condo} for $j=m+1$ implies that
the parameters $\la_{m+1}^{(a)}$ with the values $a=1,\dots,p$ can
be renumbered to satisfy the conditions
\ben
\la_{m+1}^{(a)}-\la_{m+2}^{(a)}\in\frac12\ts\ZZ_+\qquad\text{for}\qquad a=1,\dots,p.
\een
Similarly, conditions \eqref{condgl} imply that for each $i=1,\dots,m$ the parameters
$\la_{i}^{(a)}$ with $a=1,\dots,p$ can
be renumbered to satisfy the conditions
\ben
\la_{i+1}^{(a)}-\la_{i}^{(a)}\in\ZZ_+\qquad\text{for}\quad a=1,\dots,p\fand i=1,\dots,m-1.
\een
Thus, each linear weight
\ben
\la^{(a)}(u)=\big(1+\la^{(a)}_1\tss u^{-1},\dots,1+\la^{(a)}_{m+2}\tss u^{-1}\big)
\een
satisfies the conditions of Theorem~\ref{thm:necess} for $a=1,\dots,p$. By Theorem~\ref{thm:linear},
each representation $L(\la^{(a)}(u))$ is finite-dimensional and
we can conclude that $\dim L(\la(u))<\infty$ as in the proof of Corollary~\ref{cor:genhw}.

Finally, the case where some parameters $\la_{m+2}^{(a)}$ are equal to some of
$\la_m^{(k+1)},\dots,\la_m^{(p)}$ is reduced to the previous case
by considering the greatest common divisor of the polynomials
\ben
\de(u)=\gcd\big(\la_{m+2}(u),(1+\la_{m}^{(k+1)}u^{-1})\dots (1+\la_{m}^{(p)}u^{-1})\big).
\een
It suffices to observe that conditions \eqref{condo} for $j=m+1$ and \eqref{condosp}
will hold for the polynomials $\la_i(u)/\de(u)$ with $i=m,m+1$ and $m+2$ which allows us
to complete the proof of Theorem~\ref{thm:none} in this case
in the same way as above.

\medskip

\bre\label{rem:obs}
The proofs of Corollary~\ref{cor:genhw}
and Theorem~\ref{thm:none} rely on the possibility to split
the highest weight of the form \eqref{hwdec} into linear highest weights
\eqref{splinhw} satisfying the conditions of Theorem~\ref{thm:necess}.
The obstacle to prove Conjecture~\ref{conj:fdim} in this way beyond $n=1$ is the fact
that
such splitting may not be possible for the highest weights associated with
any parity sequences.

To illustrate, consider
the representation
$L(\la(u))$ of $\X(\osp_{5|2})$ associated with the standard parity sequence $\se^{\st}=100$
with the components of $\la(u)$ given by
\ben
\bal
\la_1(u)&=(1-2\tss u^{-1})(1+1/2\ts u^{-1}),\\
\la_2(u)=\la_3(u)&=(1+2\tss u^{-1})(1+1/2\ts u^{-1}),\\
\la_4(u)&=(1+3/2\tss u^{-1}).
\eal
\een
The components do satisfy the conditions of Theorem~\ref{thm:necess},
but no splitting into linear highest weights of the form \eqref{linhw}
satisfying the conditions of Theorem~\ref{thm:necess} is possible.
Moreover, such a splitting is not possible for the highest weights
corresponding to the parity sequences $\se=010$ and $001$ either.
The question whether $\dim L(\la(u))<\infty$ remains open for this example.
\ere

%
%

\bigskip
\bigskip

\small
\noindent
School of Mathematics and Statistics\newline
University of Sydney,
NSW 2006, Australia\newline
alexander.molev@sydney.edu.au

\vspace{5 mm}

\noindent
Laboratoire de Physique Th\'{e}orique LAPTh,
CNRS and Universit\'{e} de Savoie\newline
BP 110, 74941 Annecy-le-Vieux Cedex, France\newline
eric.ragoucy@lapth.cnrs.fr

\end{document}